\documentclass{article}
\usepackage{amssymb}
\usepackage{amsmath}
\usepackage{amsfonts}
\usepackage{graphicx}

\input{tcilatex}
\begin{document}

\title{On the First, the Second and the Third Stems of the Stable Homotopy
Groups of Spheres}
\author{Mehmet K\i rdar \\
Tekirda\u{g} Nam\i k Kemal University, \\
Mathematics Department, \\
Tekirda\u{g}, Turkey}
\maketitle

\begin{abstract}
We describe the first stem of the stable homotopy groups of spheres by using
some Puppe sequences, Thom complexes, K-Theory and Adams Operations
following the ideas of J. Frank Adams. We also touch upon the second and the
third stems in this perspective. Neither spectral squences nor Steenrod
operations are used.\bigskip


\end{abstract}

\section{Introduction}

The computation of the stable homotopy group of spheres $\pi _{k}^{S}=%
\underset{n\rightarrow \infty }{\lim }$ $\pi _{n+k}(S^{n})$, $k\geq 0,$ has
been one of the most interesting questions of mathematics. See [9] for a
survey. $\pi _{k}^{S}$ is called the $k$-th stem.

Skipping the zeroth stem, the purpose of this study is to make a little
survey for the first, the second and the third stems. It is well-known that
the first stem $\pi _{1}^{S}$ is $%
\mathbb{Z}
_{2}$ and it is generated by $\mathbf{\eta }$ which is the image of the
generator $\eta $ of $\pi _{1}(SO)$ under the J-homomorphism $J:\pi
_{1}(SO)\rightarrow \pi _{1}^{S}.$ The second stem $\pi _{2}^{S}$ is also $%
\mathbb{Z}
_{2}$ and it is generated by $\mathbf{\eta }^{2}$ in the ring of the stable
homotopy groups of spheres. The third stem is amazingly $%
\mathbb{Z}
_{24}$ and it is generated by\emph{\ }$\mathbf{\nu }$ which also comes from
a J-homomorphism like $\mathbf{\eta }.$

Firstly, we introduce the complex Hopf fibration and describe and compute $%
\pi _{1}^{S}$ together with $\pi _{1}(SO)$. We prove the non-triviality of
the suspension of the Hopf fibration by using complex K-theory and Adams
operations on the Puppe sequence of this fibration. Its order is also
explained by using the Puppe sequence of its double. Then, we will touch
upon the stems $\pi _{2}^{S}$ and $\pi _{3}^{S}$. We will mostly follow the
ideas in [1]. But, neither spectral sequences nor Steenrod operations will
be involved. In the appendix, we studied $\pi _{1}(SO)$ in some details.

\section{The First Stem}

One of the most interesting objects in mathematics is a Hopf fibration.
There are various definitions/constructions for them.

For example, the complex Hopf fibration $h_{c}:S^{3}\rightarrow S^{2}$ is
defined by sending a unit quaternion $q\in $ $S^{3\text{ }}$to $\overline{q}%
iq\in S^{2}\subset 
\mathbb{R}
i\oplus 
\mathbb{R}
j\oplus 
\mathbb{R}
k$ where $\overline{q}$ is the conjugate of $q.$ Fibers of this fibration
are $S^{1}$ and any two fibers are linked$.$ It can be seen as the double
cover $S^{3}\rightarrow SO(3)\rightarrow SO(3)/SO(2)=S^{2}$ of the fibration
given in the Appendix. It also can be seen as the map $S^{3}\subset 
\mathbb{C}
^{2}-\left\{ 0\right\} \rightarrow 
\mathbb{C}
P^{1}=S^{2}$ which sends $(z_{1},z_{2})$ to $\left[ z_{1},z_{2}\right] $.

Homotopy Long Exact Sequence (HLES) of this fibration gives the following

\textbf{Proposition 2.1:}

\textbf{\ i) }$\pi _{3}(S^{2})=%
\mathbb{Z}
$ generated by $\left[ h_{c}\right] $ the Hopf fibration itself.

\textbf{ii)} $\pi _{k}(S^{2})=\pi _{k}(S^{3})$ for $k\geq 3.$

The famous and deep Freudenthal suspension theorem (FST) for spheres gives
isomorphisms $\pi _{k}(S^{n})\cong \pi _{k+1}(S^{n+1})$ for $k\leq 2n-2.$
Unfortunately, $\pi _{3}(S^{2})$ is not in the stability region but $\pi
_{4}(S^{3})$ is. Therefore $\pi _{4}(S^{3})$ is the group which is
isomorphic to the first stem. On the other hand, FST for spheres gives
epimorphisms $\pi _{2n-1}(S^{n})\rightarrow \pi _{2n}(S^{n+1})$ so that we
have an epimorphism $\pi _{3}(S^{2})\rightarrow \pi _{4}(S^{3}).$ Therefore, 
$\pi _{4}(S^{3})$ is either trivial or a nontrivial cyclic group generated
by $\left[ Sh_{c}\right] $. In fact, according to the homotopy table of [7]
on page 339, it is nontrivial and $\pi _{4}(S^{3})=%
\mathbb{Z}
_{2}$. Therefore, the latter argument is true and we are done.

But, we want to explain the non-triviality of $\left[ Sh_{c}\right] $ and
the relation $2\left[ Sh_{c}\right] =0$ without using the information $\pi
_{4}(S^{3})=%
\mathbb{Z}
_{2}$.

\textbf{Proposition 2.2: }$\left[ Sh_{c}\right] \neq 0.$

\textit{Proof: }The suspension of the mapping cone of the map $h_{c}$, that
is $S(S^{2}\cup _{h_{c}}B^{4}),$ \textit{\ }is homotopy equivalent to\textit{%
\ }the mapping cone of the map $Sh_{c}$, that is to $S^{3}\cup
_{Sh_{c}}B^{5},$ due to the Puppe sequence of $h_{c}$, [1]. But, it is
well-known that the mapping cone $S^{2}\cup _{h_{c}}B^{4}$ is the Thom
complex of the Hopf line bundle $\eta $ over $%
\mathbb{C}
P^{1}$ and it is the complex projective space $%
\mathbb{C}
P^{2}$. See [3] for Thom complexes. If $Sh_{c}$ was trivial, $S%
\mathbb{C}
P^{2}$ would be split, that is, it would be the wedge sum $S^{3}\vee S^{5}$
and thus $S^{2}%
\mathbb{C}
P^{2}$ would be $S^{4}\vee S^{6},$ etc.

We must show that $S%
\mathbb{C}
P^{2}$ and $S^{3}\vee S^{5}$ or $S^{2}%
\mathbb{C}
P^{2}$ and $S^{4}\vee S^{6}$ are distinct, that is, not homotopy equivalent.
The reader may think that these distinctions are obvious. In fact, for
example, they are a very special case of Theorem 1.1. of [6], S. Feder and
S. Gitler's article about the classification of the stable homotopy types of
the stunted projective spaces.

But I will try to do an awkward proof for the latter distinction.

By Bott periodicity, we have $\widetilde{K}(S^{2}%
\mathbb{C}
P^{2})=\widetilde{K}(%
\mathbb{C}
P^{2})\otimes \widetilde{K}(S^{2})=\dfrac{%
\mathbb{Z}
(\mu \nu )\oplus 
\mathbb{Z}
(\mu ^{2}\nu )}{\left( \nu ^{2},\mu ^{3}\right) }$ where $\mu $ is the
reduction of the complex Hopf line bundle $\eta _{c}$ over $%
\mathbb{C}
P^{2}$ and $\nu \in \widetilde{K}(S^{2})$ is the Bott element .On the other
hand, $\widetilde{K}(S^{4}\vee S^{6})=\widetilde{K}(S^{4})\oplus \widetilde{K%
}(S^{6})=%
\mathbb{Z}
(x)\oplus 
\mathbb{Z}
(y)$ where virtually $\dim x=4,$ $\dim y=6.$ These rings are really
isomorphic.

Suppose we an homotopy equivalence $\phi :S^{2}%
\mathbb{C}
P^{2}\rightarrow S^{4}\vee S^{6}$ so that we have an induced ring
isomorphism $\phi ^{\ast }:$ $\widetilde{K}(S^{4})\oplus \widetilde{K}%
(S^{6})\rightarrow $ $\widetilde{K}(S^{2}%
\mathbb{C}
P^{2}).$ We may assume that $\phi ^{\ast }(x)=+\mu \nu .$

Now, we must have $\phi ^{\ast }(\psi ^{p}(x))=\psi ^{p}(\phi ^{\ast }(x))$
due the neutrality of the complex Adams operation $\psi ^{p}$ for some
natural number $p>1.$ But then, by properties of Adams operations, [6], $%
\phi ^{\ast }(p^{2}x)=\psi ^{p}(\mu \nu )$ and $p^{2}\mu \nu =\psi ^{p}(\mu
)\psi ^{p}(\nu )=p\nu \psi ^{p}(\mu ).$ Hence we get $\psi ^{p}(\mu )=p\mu .$
But this is a contradiction. Because, $\psi ^{p}(\mu )=p\mu +\frac{p(p-1)}{2}%
\mu ^{2}.$

Next Proposition is also well-known to experts.

\textbf{Proposition 2.3: }$2\left[ Sh_{c}\right] =0.$

\textit{Proof: }$S^{2}\cup _{2h_{c}}B^{4}$ is the Thom space of $\eta ^{2}$
like $S^{2}\cup _{h_{c}}B^{4}$ is the Thom space of $\eta .$

But, in $KO(S^{2})=KO(%
\mathbb{C}
P^{1})$ we calculate $r(\eta ^{2})+2=2r(\eta )=4$ and it is trivial.

Therefore, it follows that $S^{N}T(r(\eta ^{2})+2)=S^{N}T(4)$ for some
natural number $N$ where equality means the homotopy equivalence$.$

Hence we get, $S^{N+2}(S^{2}\cup _{2h_{c}}B^{4})=S^{N+4}\vee S^{N+6}$. 

On the other hand, due to Puppe sequence of $2h_{c},$ we have $%
S^{N+2}(S^{2}\cup _{2h_{c}}B^{4})=S^{N+4}\cup _{2S^{N+2}h_{c}}B^{N+6}.$ 

Thus $2\left[ Sh_{c}\right] =2\left[ S^{N+2}h_{c}\right] =\left[
2S^{N+2}h_{c}\right] =0$.

We must remark here that what makes the splitting possible is actually the $%
J $-order of the bundle, not the $KO$-order, in the above proposition. But,
since these orders are equal, we have no problem here. But, in the third
stem, it will be different. We also note that we did not mention about the $%
J $-homomorphism so far!

\textbf{Corollary 2.4: }$\pi _{4}(S^{3})=%
\mathbb{Z}
_{2}$ and it is generated by $\left[ Sh_{c}\right] .$

\textbf{Corollary 2.5:} $\pi _{1}^{S}=%
\mathbb{Z}
_{2}$ and it is generated by $\mathbf{\eta }$, the infinite suspension of
the Hopf fibration.

J- Homomorphism proof of Corollary 2.5 is briefly explained in the following
way: From a map $S^{r}\longrightarrow O(n),$ $r\geq 0,$ one can construct a
map $S^{r+n}\rightarrow S^{n}$ by the "Hopf construction". One can show that
this construction results in a homomorphism $J:\pi _{r}(O(n))\rightarrow \pi
_{r+n}(S^{n})$ and also by taking limit, a homomorphism $J:\pi _{r}(O)=%
\widetilde{KO}(S^{r+1})\rightarrow \pi _{r}^{S}.$ And most interestingly,
this homomorphism factors through a group denoted by $\widetilde{J}%
(S^{r+1}), $ the reduced $J$-group of the sphere $S^{r+1}.$

My supervisor \.{I}brahim Diba\u{g} devoted great efforts for determining
the decomposition of $\widetilde{J}$-groups of complex projective spaces and
lens spaces as a direct-sum of cyclic groups. These groups involve very
complicated arithmetic. I think his efforts culminated in [5].

In our case $r=1$, one has the homomorphism $J:\pi _{1}(O)=\pi
_{1}(SO)\rightarrow \pi _{1}^{S}.$ And furthermore, one can prove, that $%
J(\eta )=\mathbf{\eta }$ where $\eta $ is the generator of $\pi _{1}(SO)$ as
explained in the Appendix. Hence, our isomorphism is achieved.

Adams, [1], introduced the $d$ and $e$ invariants to detect the image of $\ $%
the $J$-homomorphism for all $r$. He proved that for $r>0$ and $r\equiv 0,$ $%
1$ $(\func{mod}$ $8)$, \ $J$-homomorphism is injective and gives a $%
\mathbb{Z}
_{2}$ direct sum part of the $r$-th stem.

\section{The Second Stem}

FST implies that $\pi _{2}^{S}=\pi _{6}(S^{4}).$ But according to the
homotopy table of [7] again, $\pi _{6}(S^{4})=%
\mathbb{Z}
_{2}$. Therefore, easily $\pi _{2}^{S}=%
\mathbb{Z}
_{2}.$

On the other hand, this stem is a phenomenon of the 20th century. First of
all, $\widetilde{KO}(S^{3})=0$ and the generator of the second stem is not
image of the $J$-homomorphism. But, this generator is the direct limit of
the composite $S^{n}h_{c}\circ S^{n+1}h_{c}$ when $n\longrightarrow \infty .$
In other words, it is $\mathbf{\eta }^{2}$ in the ring of stable homotopy
group of spheres $\dbigoplus\limits_{k=0}^{\infty }\pi _{k}^{S}$. But, we
know that $\mathbf{\eta }$ is a image of the $J$-homomorphism. This is weird.

By the epimophism part of FST, it is enough to understand the composite $%
Sh_{c}\circ S^{2}h_{c}.$

We can explain the order: $2[Sh_{c}\circ S^{2}h_{c}]=[2Sh_{c}\circ
S^{2}h_{c}]=[0\circ S^{2}h_{c}]=0.$

But, the main problem of this stem is to show the nontriviality of this
composite.

\textbf{Proposition 3.1: }$[Sh_{c}\circ S^{2}h_{c}]\neq 0.$

\textit{Proof: }HLES of the complex Hopf fibration shows that $\pi _{\ast
}:\pi _{4}(S^{3})\rightarrow \pi _{4}(S^{2})=%
\mathbb{Z}
_{2}$ is isomorphism and since $\pi _{4}(S^{3})$ is generated by $\left[
Sh_{c}\right] ,$ we have that $\pi _{4}(S^{2})$ is generated by $\left[
h_{c}\circ Sh_{c}\right] .$ On the other hand, we have a 2-local fibration $%
S^{2}\rightarrow \Omega S^{3}\rightarrow \Omega S^{5}$ which gives the
so-called EHP-sequence, [9], and thus an isomorphism $E:\pi
_{4}(S^{2})\rightarrow \pi _{5}(S^{3})$ where $E$ is $S,$ the suspension
homomorphism. The result follows since $E(\left[ h_{c}\circ Sh_{c}\right] )=%
\left[ Sh_{c}\circ S^{2}h_{c}\right] $.

The proof above is quite unsatisfactory for me. Because, I want to try to
find an explanation by means of the mapping cone space $M=S^{3}\cup
_{f}B^{6} $ and related maps with it, where $f=Sh_{c}\circ S^{2}h_{c}.$ Is $%
M $ the 6-th skeleton $\left( \Omega S^{5}\right) ^{\left( 6\right) }$
2-locally? The attaching map $f$ and the space $M$ are mysterious.

\textbf{Corollary 3.2: }$\pi _{2}^{S}=%
\mathbb{Z}
_{2}$ and it is generated by $\mathbf{\eta }^{2}.$

An alternative proof of Proposition 3.1 can be obtained from Theorem 7.2 of
[1]. Adams claims that his $d_{%
\mathbb{R}
}$ invariant detects the nontriviality of $\mathbf{\eta }^{2}.$

Another explanation for the non-triviality of $\mathbf{\eta }^{2}$ is the
Pontrjagin approach, [9]: $\pi _{k}^{S}$ is isomorphic to the cobordism
group of framed $k$-dimensional manifolds in $%
\mathbb{R}
^{\infty },$ denoted by $\Omega _{k}^{fr}$. Under the isomorphism $\pi
_{k}^{S}$ $\cong \Omega _{k}^{fr}$, $\mathbf{\eta }$ corresponds to the
class of the circle $\left[ S^{1}\right] $ and $\mathbf{\eta }^{2}$
corresponds to the class of the torus $\left[ S^{1}\times S^{1}\right] .$
The non-triviality of $\mathbf{\eta }^{2}$ is detected by an invariant
called the Arf-Kervaire invariant which is the Arf invariant of the
quadratic form on the mod 2 first homology of the torus as far as I
understand. You may see the famous Turkish mathematician Cahit Arf's paper
for the origin of this invariant, [2].

\section{The Third Stem}

We have the following ingredients: With the help of homotopy tables, FST
implies $\pi _{3}^{S}=\pi _{8}(S^{5})=%
\mathbb{Z}
_{24}$. There is a quaternionic Hopf Fibration $h_{q}:S^{7}\rightarrow S^{4}$
with fiber $S^{3}$ and its suspension is the generator of $\pi _{8}(S^{5}),$
so its infinite suspension $\mathbf{\nu }$ is the generator of $\pi
_{3}^{S}. $ Furthermore, it is the image of $x=1\in \widetilde{KO}(S^{4})=%
\mathbb{Z}
$ under the $J$-homomorphism. Notice that $x$ has infinite order while $%
\mathbf{\nu }$ has finite order.

Why is $d\left[ Sh_{q}\right] =0$ (or $d\mathbf{\nu }=0)$ if and only if $24$
divides $d?$ We can mimic the proof of Proposition 2.2 to deduce that $\left[
Sh_{q}\right] \neq 0$ in the following way: We have a quaternionic Hopf line
bundle $\eta _{q}$ over the quaternionic projective space $\mathbf{H}%
P^{1}=S^{4}.$ The mapping cone of $h_{q}$ is the Thom complex $T(\eta _{q})=%
\mathbf{H}P^{2},$ the quaternionic projective plane. Then, by looking at $%
\widetilde{K}(S^{4}\mathbf{H}P^{2})$ and the effect of Adams operations on
it, we may deduce that $S^{5}\cup _{Sh_{q}}B^{9}$ is not homotopy equivalent
to $S^{5}\vee S^{9}.$ But this solves only the case $d=1$ of the question we
just stated.

Instead, we will use the Adams conjecture to find the upper bound. It is
such a deep theorem that it is still called as a conjecture.

\textbf{Proposition 4.1: }The order of $\mathbf{\nu }$ divides $24.$

\textit{Proof: }You may see Diba\u{g}'s article [4] for a short elementary
and self-contained proof of the Adams conjecture. For $CW$-complexes, it
implies for sufficiently big natural number $N,$ 
\begin{equation*}
k^{N}(\psi _{%
\mathbb{R}
}^{k}(x)-x)=0
\end{equation*}%
in the $\widetilde{J}$-group $\widetilde{J}(X)$ for all $x\in \widetilde{KO}%
(X)$ and natural numbers $k$.

If we take $X=S^{4},$ $x=1\in \widetilde{KO}(S^{4})=%
\mathbb{Z}
,$ we get $k^{N}(k^{2}-1)x=0$ in $\widetilde{J}(S^{4})$ $.$ The common
multiple of the numbers $k^{N}(k-1)(k+1)$ is 24 while $N$ is big. So, $24x=0$
in $\widetilde{J}(S^{4})$ and under the $J$-homomorphism we get $24\mathbf{%
\nu }=0$.

The Proposition above implies that the space $S^{5}\cup _{24Sh_{q}}B^{9}$ is
split, hence we get $S(S^{4}\cup _{24h_{q}}B^{8})=S^{5}\vee S^{9}.$On the
other hand, $S^{4}\cup _{24h_{q}}B^{8}$ is the Thom space $T(\eta
_{q}^{24}). $ In $KO(S^{4})=KO(\mathbf{H}P^{1}),$ we have $\eta
_{q}^{24}+92=24\eta _{q}. $ The Thom space $T(24\eta _{q})$ is the stunted
quaternionic projective space $\mathbf{H}P^{25}\diagup \mathbf{H}P^{23},$ [3]%
$.$ It follows from these facts that for some natural number $N$ we have 
\begin{equation*}
S^{N}(\mathbf{H}P^{25}\diagup \mathbf{H}P^{23})=S^{96+N}\vee S^{100+N}
\end{equation*}%
as an homotopy equivalence equality.

For the lower bound, we must prove the following Proposition. We will refer
to [6] mentioned in Proposition 2.2.

\textbf{Proposition 4.2: }$12\mathbf{\nu \neq 0.}$

\textit{Proof: }By a similar argument made just after the Preposition 4.1,
we deduce that if $12\mathbf{\nu }$ were zero then $\mathbf{H}P^{13}\diagup 
\mathbf{H}P^{11}$ would be $S^{48}\vee S^{52}.$ We must show that this is
not possible. This can easily follow as a very special case of Theorem 1.2
of [6]: The spaces $\mathbf{H}P^{n+k}\diagup \mathbf{H}P^{k-1}$ and $\mathbf{%
H}P^{n+l}\diagup \mathbf{H}P^{l-1}$ have the same stable homotopy type if $%
k-l=O(B_{n})$ where $B_{n}$ is the $J$-order of the Hopf bundle of $\mathbf{H%
}P^{n}.$ If we take $n=1,$ $B_{1}=24.$ Now, let $k=12$ and $l=0.$ Then it
follows that $\mathbf{H}P^{13}\diagup \mathbf{H}P^{11}$ can not be stably
homotopic to $\mathbf{H}P_{+}^{1}=S^{0}\vee S^{4}$ which is stably homotopic
to $S^{48}\vee S^{52}.$

\textbf{Corollary 4.3: }$\pi _{2}^{S}=%
\mathbb{Z}
_{24}$ and it is generated by $\mathbf{\nu }.$

Adams, [1], generalizes the last result as $\widetilde{J}(S^{4k})=%
\mathbb{Z}
_{m(2k)}$ where $m(2k)$ is the denominator of $\dfrac{B_{2k}}{4k}$ and $%
B_{2k}$ is the $2k$-th Bernoulli number and this group is a direct summand
of $\pi _{4k-1}^{S}.$

Pontrjagin approach is mind blowing at this stem, [9]. Under the isomorphism 
$\pi _{3}^{S}$ $\cong \Omega _{3}^{fr}$, $\mathbf{\nu }$ corresponds to the
class of the three dimensional sphere $\left[ S^{3}\right] $ and it is told
that deletion of $24$ three dimensional balls around the singularities of a $%
K3$-surface corresponds to the null cobordism $24\left[ S^{3}\right] =0.$

As a final note, I recall some relations. Two of them are $\eta ^{4}=\eta
\nu =0.$ If we cheat a little without giving a constructive proof, these
relations follow from the fact that $\pi _{4}^{S}=0.$ Another relation is $%
\eta ^{3}=12\nu .$ The explanation of this relation without spectral
sequences must be hard.

\section{Appendix: The Group $\protect\pi _{1}(SO)$}

$\allowbreak $Let $SO(n)$ be the group of $n\times n$ special orthogonal
matrices. For each $n\in 
\mathbb{Z}
^{+},$ we have natural group inclusions $i_{n}:SO(n)\rightarrow SO(n+1)$
which give fibrations $SO(n)\rightarrow SO(n+1)\rightarrow S^{n}$ by group
quotient. Hence, due to HLES of these fibrations, the maps $i_{n}$ are $%
(n-1) $-connected maps so that $\pi _{k}(SO(n))=\pi _{k}(SO)$ for all $n\geq
k+2$ where $SO$ is the infinite dimensional special orthogonal group,
obtained by taking the direct limit of the group $SO(n)$ when $n\rightarrow
\infty $. Similar arguments goes for the orthogonal group $O(n)=%
\mathbb{Z}
_{2}\times SO(n)$ and the infinite dimensional orthogonal group $O=%
\mathbb{Z}
_{2}\times SO$. Note that $\pi _{k}(O)=\pi _{k}(SO)$ for all $k\geq 1$ but $%
\pi _{0}(SO)=0$ whereas $\pi _{0}(O)=%
\mathbb{Z}
_{2}.$

For $k=1,$ we must take at least $n=3$ to calculate $\pi _{1}(SO)=%
\mathbb{Z}
_{2}.$ Really, $\pi _{1}(SO(1))=0$ and $\pi _{1}(SO(2))=%
\mathbb{Z}
$, therefore this bound is strict.

\textbf{Proposition A.1: }$SO(3)=%
\mathbb{R}
P^{3}.$

\textit{Proof:} $S^{3}\subset \mathbf{H}$ can be considered as the space of
unit quaternions and each $q\in $ $S^{3}$ defines a rotation $x\rightarrow 
\overline{q}xq$ about the origin of the real vector space of pure imaginary
quaternions $%
\mathbb{R}
i\oplus 
\mathbb{R}
j\oplus 
\mathbb{R}
k\cong 
\mathbb{R}
^{3}$. This map turns out to be a surjective group homomorphism $%
S^{3}\rightarrow SO(3)$ with kernel $\left\{ \pm 1\right\} .$Hence we get an
isomorphism $S^{3}/%
\mathbb{Z}
_{2}\rightarrow SO(3).$

\textbf{Corollary A. 2: }$\pi _{1}(SO)=%
\mathbb{Z}
_{2}.$

\textit{Proof:} HLES of the fibration $%
\mathbb{Z}
_{2}\rightarrow S^{3}\rightarrow SO(3)$ gives $\pi _{1}(SO(3))=\pi _{0}(%
\mathbb{Z}
_{2})=%
\mathbb{Z}
_{2}.$The result follows by stability.

Another but a more difficult proof of Proposition A.1 can be given in the
following way. Each rotation of $%
\mathbb{R}
^{3}$ about the origin is determined by a vector which is called the axis of
the rotation and a rotation angle $\theta \in \lbrack 0,\pi ].$ When the
rotation angle is $\pi ,$ antipodal directions give the same rotation.
Consider the closed ball in $%
\mathbb{R}
^{3}$ with center the origin and radius $\pi .$ Identify the antipodal
points of the boundary and get a quotient space $X.$ $X$ is homeomorphic to $%
\mathbb{R}
P^{3}.$ Define $\phi :X\rightarrow SO(3)$ by sending $[x,y,z]\in X$ to the
rotation with direction determined by the vector $(x,y,z)$ and rotation
angle $\sqrt{x^{2}+y^{2}+z^{2}}.$ This map is a homeomorphism.

The generator of the group $\pi _{1}(SO(3))$ is described by using the space 
$X$ and the map $\phi $ above in the following way. Consider the loop $%
\gamma :I=\left[ 0,1\right] \rightarrow X$ defined by $\gamma (t)=\left[
0,0,\pi \cos \pi t\right] $.

\textbf{Proposition A.3: }$\phi _{\ast }(\left[ \gamma \right] )$ is the
generator of $\pi _{1}(SO(3)).$

Proof of \ Proposition A.3 will be given later in this section. The loop $%
\gamma $ corresponds to the set of rotations about $z$-axis, in other words,
image of $\gamma $ in $SO(3)$ is exactly $SO(2).$ The following proposition
is equivalent to the Proposition A.3 and it is an alternative proof.

\textbf{Proposition A.4: }$\pi _{1}(SO(3))=<i_{2_{\ast }}(1)>$ where $1\in
\pi _{1}(SO(2))=%
\mathbb{Z}
$.

\textit{Proof:} HLES of the fibration $SO(2)\rightarrow SO(3)\rightarrow
S^{2}$ gives the epimorphism $\pi _{1}(SO(2))=%
\mathbb{Z}
\rightarrow \pi _{1}(SO(3))=%
\mathbb{Z}
_{2}\rightarrow 0.$

In other words, if we identify $SO(2)$ with $S^{1},$ the natural inclusion $%
i_{2}:SO(2)\subset SO(3)$ and similarly $i:SO(2)\subset SO$ can be seen as
the generators of the corresponding fundamental groups. Let $\eta =\left[ i%
\right] \in \pi _{1}(SO).$ Then $\eta $ is the generator of $\pi _{1}(SO)$
and it follows from Corollary A.2 that

\textbf{Corollary A.5: }$2\eta =0.$

But, the geometric explanation of Corollary A.5 is quite complicated. It is
enough to show that $2[\gamma ]=0$ in $\pi _{1}(X).$

Define the loops $\alpha ,\beta :I\rightarrow X$ by $\alpha (t)=\left[
0,-\pi \sin \pi t,\pi \cos \pi t\right] $ and $\beta (t)=\left[ 0,\pi \sin
\pi t,\pi \cos \pi t\right] .$ Then $\left[ \beta \right] =-\left[ \alpha %
\right] .$ On the other hand $\left[ \gamma \right] =\left[ \alpha \right] $
by the ellipse shaped homotopy $H(s,t)=\left[ 0,-\pi s\sin \pi t,\pi \cos
\pi t\right] $ and $\left[ \gamma \right] =\left[ \beta \right] $ by a
similar homotopy $H(s,t)=\left[ 0,\pi s\sin \pi t,\pi \cos \pi t\right] .$
Hence the result follows.

We can write special orthogonal matrices given their unit directions and
angles. If we do that, the paths $\alpha ,\beta $ correspond to paths 
\begin{equation*}
\left[ 
\begin{array}{ccc}
-1 & 0 & 0 \\ 
0 & -\cos 2\pi t & -\sin 2\pi t \\ 
0 & -\sin 2\pi t & \cos 2\pi t%
\end{array}%
\right] \text{ and }\left[ 
\begin{array}{ccc}
-1 & 0 & 0 \\ 
0 & -\cos 2\pi t & \sin 2\pi t \\ 
0 & \sin 2\pi t & \cos 2\pi t%
\end{array}%
\right]
\end{equation*}%
in $SO(3),$ respectively. They are inverse paths to each other as expected.
The path $\gamma $ of course corresponds to $\left[ 
\begin{array}{ccc}
\cos \pi t & \sin \pi t & 0 \\ 
-\sin \pi t & \cos \pi t & 0 \\ 
0 & 0 & 1%
\end{array}%
\right] $. One can also write the homotopy matrices from $\gamma $ to $%
\alpha $ or $\beta $ in $SO(3).$

\textit{Proof of Proposition A.3}: It is enough to show that $\left[ \alpha %
\right] $ is not trivial. Suppose that it is trivial with an homotopy to the
constant loop at the pole. By compressing this homotopy to the boundary we
obtain a homotopy in the boundary. But the boundary is $%
\mathbb{R}
P^{2}$ and thus $\left[ \alpha \right] $ is trivial in $%
\mathbb{R}
P^{2}.$ Now, a similar argument can be repeated for $%
\mathbb{R}
P^{2}$concluding that $\left[ \alpha \right] $ is trivial in $%
\mathbb{R}
P^{1}=S^{1}$. But, $\left[ \alpha \right] $ is the generator of the
fundamental group of this last boundary by construction. A contradiction.

Vector bundle interpretation can be stated as $\pi _{1}(SO)=\pi _{1}(O)=\pi
_{2}(BO)=\widetilde{KO}(S^{2})=%
\mathbb{Z}
_{2}.$ There is a realification homomorphism $r:\widetilde{K}%
(S^{2})\rightarrow \widetilde{KO}(S^{2}).$ Since $S^{2}\cong 
\mathbb{C}
P^{1},$ we have $\widetilde{K}(S^{2})=\widetilde{K}(%
\mathbb{C}
P^{1})=%
\mathbb{Z}
$ generated by $\eta -1,$ reduction of the Hopf line bundle $\eta :E(\eta
)\rightarrow 
\mathbb{C}
P^{1}.$ The homomorphism $r$ is surjective and hence $\widetilde{KO}(S^{2})$
is generated by $r(\eta )-2,$ the reduction of the realification of the Hopf
line bundle and its order in this group is 2.

\bigskip

\end{document}